\documentclass[12pt]{amsart}

\usepackage{amsmath,amssymb,amscd,amsfonts}

\newtheorem{thm}{Theorem}[section]
\newtheorem{lem}[thm]{Lemma}

\newtheorem{prop}[thm]{Proposition}

\newtheorem{assu-nota}[thm]{Assumption--Notation}

\theoremstyle{remark}
\newtheorem{ex}{Example}
\newtheorem{remark}{Remark}

\newcommand{\C}{\mathbb C}
\newcommand{\Z}{\mathbb Z}
\newcommand{\Q}{\mathbb Q}
\newcommand{\R}{\mathbb R}
\newcommand{\pp}{\mathbb P}
\newcommand{\Aut}{\operatorname{Aut}}

\newcommand{\Pic}{\operatorname{Pic}}
\newcommand{\Num}{\operatorname{Num}}
\newcommand{\OO}{\mathcal O}

\newcommand{\epsi}{\epsilon}

\newcommand{\Ga}{\Gamma}
\newcommand{\Si}{\Sigma}
\newcommand{\fie}{\varphi}

\newcommand{\inv}{^{-1}}

\numberwithin{equation}{section}

\newcommand{\Hom}{\text{Hom}}

\newcommand{\calC}{\mathcal{C}}

\newcommand{\calO}{\mathcal{O}}

\newcommand{\bbF}{\mathbb{F}}
\newcommand{\bbC}{\mathbb{C}}

\newcommand{\bbP}{\mathbb{P}}
\newcommand{\bbQ}{\mathbb{Q}}
\newcommand{\bbR}{\mathbb{R}}

\newcommand{\bbZ}{\mathbb{Z}}
\newcommand{\bbG}{\mathbb{G}}

\newcommand{\bfF}{\mathbf{F}}
\newcommand{\F}{\mathbf{F}}

\begin{document}
\title{Rational surfaces with many nodes}
\author{Igor Dolgachev, Margarida Mendes Lopes, Rita Pardini} 
\thanks{Research of the first author is partially supported by a NSF grant, the second author
is a member of CMAF and of the Departamento de Matem\'atica da Faculdade
de Ci\^encias da Universidade de Lisboa and the third author is a member
of GNSAGA of CNR}
\maketitle
\section{Introduction}
Let $X$ be a smooth rational projective algebraic surface over an
algebraically closed
field ${\bf k}$ of characteristic $\ne 2$. It is
known that for any nodal curve
$C$ on $X$  there
exists a birational morphism $f:X\to X'$ such that the image of $C$ is an
ordinary double point (a {\em node}). Let $n(X)$ be the maximal number of
disjoint nodal curves on $X$. After blowing down all of them we obtain
a rational normal surface $X'$ with $n(X)$ nodes. The Picard number
$\rho(X')$ of
$X'$ is equal to the Picard number $\rho(X)$ of $X$ minus $n(X)$. Since
$X'$ is projective, $\rho(X') = \rho(X)-n(X) \ge 1$. In this paper we
study the limit cases, namely $\rho(X') = 1$ or $2$. More precisely, we
prove that
$\rho(X') = 1$ is possible only  if $X'$ is isomorphic to a
quadric cone and we describe all the $X'$'s   such
that
$\rho(X') = 2$.

The question of the number of nodes on an algebraic surface is a very old
one
and   has a long history, but, to our knowledge,  this particular
problem has not been considered. Our interest in this question arose in the
course of
investigating complex surfaces of general type with $p_g=0$ admitting a
double
plane construction, and  in the last section of this
paper, working over $\bbC$, we give an application to such surfaces with
$K^2=8,9$. More precisely, we extend some of the
results of the previous sections for surfaces with
$p_g=q=0$ and  non-negative Kodaira dimension and then we
consider surfaces
$S$ of general type with $p_g=0$ with an involution
$\sigma$. We show that if
$K^2=9$, then $S$ does not admit an involution $\sigma$ and
we list all the possibilities for the quotient surface
$S/<\sigma>$\, if
$K^2=8$.

One of our  main tools is the code associated to a set of nodal curves (see
Section~\ref{codes}), which
has already been considered by A. Beauville in \cite{beauville}.

\medskip
\paragraph{\bf Notations and conventions} As already explained,  we work
over any
algebraically closed field
${\bf k}$ of characteristic
$\ne 2$ in sections $2$ and $3$, whilst  in   section $4$ we work over
$\bbC$.

The multiplicative group of ${\bf k}$ is denoted by $\bbG_m$. For any abelian
group $A$ we denote by $_2A$ the kernel of the homomorphism
$[2]\colon A\to A$, $a\mapsto 2a$.

 All varieties are projective algebraic. We do not
distinguish between line bundles and divisors on a smooth
variety, and use additive and multiplicative notation
interchangeably. Linear equivalence is denoted by $\equiv$
and numerical equivalence by $\sim$. The intersection
product of divisors (line bundles) $A$ and $B$ is denoted  by
$A\cdot B$. We denote by $\kappa(X)$ the Kodaira dimension of a variety $X$
and by $\rho(X)$ the Picard number of  $X$, that
is the rank of
the N\'eron--Severi group of $X$. A {\it nodal} curve on a
surface is a   smooth rational
curve $C$  such that $C^2=-2$.
 The remaining notation is standard in
algebraic geometry.

\section{Nodal curves, binary codes and covers}\label{codes}

In this section all varieties are defined over an algebraically closed field ${\bf k}$
of characteristic $\ne 2$.

 Recall that a binary code is a subspace $V$ of a
$k-$dimensional vector space $W$ over $\bbF_2$ equipped with a basis
$(e_1,\ldots,e_k)$. The dimension of $W$ (identified in the sequel with
$\bbF_2^k$) is called the \emph{length} of $V$. For each
$v\in V$ the number of nonzero coordinates of
$v$ with respect to the basis is called the \emph{weight} of $v$.\par
Two codes $V_1, V_2\subset \bbF_2^k$ are \emph{isomorphic} if
there exists a permutation of the coordinates of $\bbF_2^k$ mapping $V_1$
onto $V_2$.

We say that a code $V$ is {\em reduced} if there is no $1\le i\le k$ such
that
$V\subset \{x_i=0\}$. To every code $V$ one can associate a reduced code
$V'$, simply by
deleting the ``useless'' coordinates. The  dimension and  the weights of $V$
and
$V'$ are
the same, while $V'$ has (possibly) smaller length. We say that two codes
$V_1\subset
\bbF_2^{k_1}$, $V_2\subset
\bbF_2^{k_2}$ are {\em essentially isomorphic} if the corresponding reduced
codes are
isomorphic.

We mention here a code that plays an important part in what follows. Given
an integer
$n$,   consider the code of even vectors $V=\{\sum x_i=0\}\subset
\bbF_2^n$.
$V$ has
dimension
$n-1$ and its weights are all even. We define the {\em code of   doubly even
vectors}
$DE(n)$ to be the image of $V$ via the injection $\bbF_2^n\to\bbF_2^{2n}$
defined by
$(x_1\ldots x_n)\mapsto (x_1x_1\ldots x_nx_n)$. So $DE(n)$ has length $2n$,
dimension
$n-1$ and all its weights are divisible by $4$.\par
It is possible to associate to a linear code $V\subset \bbF_2^n$  a
lattice
$\Gamma_V$ in the Euclidean space
$\bbR^n$ (see, for example,
\cite{ebeling}).  One considers the canonical homomorphism $p:\bbZ^n\to
\bbF_2^n$ and takes $\Gamma_V$ to be $\frac{1}{\sqrt{2}}p^{-1}(V)$.  For
example, the code
$V$ of even vectors in $\bbF_2^n$
defines the root lattice of type $D_n$. The code of doubly even
vectors $DE(n)$ defines the root lattice $D_{2n}$ (loc. cit., p.25). \par

 Binary
codes arise naturally in the theory of algebraic surfaces, as follows.
Consider a smooth projective surface $Y$
 and $k$
 disjoint  nodal curves  $C_1,\ldots C_k$ of $Y$. Let $\calC$ be
the subgroup of
$\Pic(Y)$ generated by the curves $C_i$, which  is a free abelian group of
rank $k$.  Let
$\fie:\calC/2\calC\to
\Pic(Y)/2\Pic(Y)$ be the natural homomorphism of $2-$elementary abelian
groups. We call the kernel $V$ of $\fie$ the
{\em (binary) code associated to the $C_i$} and denote its dimension by
$r$. Here we take for a basis of
$W := \calC/2\calC$ the classes of the curves $C_i$ modulo $2\calC$.
The
lattice $\Gamma_V$ is
isomorphic to the smallest primitive sublattice containing the $C_i$ of the
lattice
$\Num(X)$ of divisors of $X$  modulo numerical equivalence.

We say that a curve $C_i$ {\em appears in $V$}  if $V$ is not contained in
$\{x_i=0\}$ and we denote by $m$ the number of $C_i$ that appear in $V$ (so
$m$
is the length of the reduced code associated to $V$).  The vector
$v=(x^1\ldots x^k)\in\bbF_2^k$ is in
$V$ if and only if  there exists $L_v\in \Pic(Y)$ such that
$2L_v\equiv \sum x^iC_i$ (when it is convenient, we identify   $0, 1\in
\bbF_2$
with the integers
$0,1$). Notice
that $K_Y\cdot L_v=0$ and thus $L_v^2$ is even by the adjunction formula.
Then the
weight $w(v)$ of $v$ is equal to
$-2L_v^2$ and so it is divisible by $4$.
Notice   that
$L_v$ is uniquely determined by $v$ if and only if  $_2\Pic(Y)=0$.

The
following result is analogous to the construction of the Galois cover  of a
surface $Y$ associated to a torsion subgroup of
$\Pic(Y)$.

\begin{prop}\label{cover} Let $Y$ be a smooth projective surface with\,
$_2\Pic(Y)=0$, let $C_1\ldots C_k$ be disjoint nodal curves
of
$Y$ and let $V$, $L_v$ be defined as above.

Then there exists a unique
smooth connected Galois cover
$\pi\colon Z\to Y$ such that:
\begin{itemize}
\item[(i)]  the Galois group of $\pi$ is $G:=\Hom(V,\bbG_m)$;

\item[(ii)]  the branch locus of $\pi$ is the union of the $C_i$ that appear in $V$;

\item[(iii)]  $\pi_*\OO_Z=\oplus_{v\in V}L_v\inv$, and $G$ acts on $L_v\inv$ via the
character
$v\in V\cong Hom(G,\bbG_m)$.
\end{itemize}
\end{prop}
\begin{proof} For $v\in V$ and $g\in G$, we define $\epsi_v(g)\in\{0,1\}$ by
$(-1)^{\epsi_v(g)}=g(v)$.  We fix a basis $v_1\ldots v_r$ of
$V$ and we write $\epsi_j$ for $\epsi_{v_j}$, $j=1\ldots r$.

By  Proposition $2.1$ of
\cite{ritaabel},   in order to determine
$\pi\colon Z\to Y$ we have to assign the (reduced) {\em building data},
namely:

1) for every
nonzero $g\in G=\Hom(V,\bbG_m)$ an effective  divisor $D_g$;

2) for every $j=1\ldots r$ a line bundle
$M_j$

\noindent in such a way that the following relations are satisfied:
\begin{equation}\label{fund}
2M_j\equiv \sum_{g\in G} \epsi_j(g)D_g, \quad j=1\ldots r.
\end{equation}
For $i=1\ldots k$ we denote by $\psi_i\colon W=\bbF_2^k\to \bbG_m$ the
homomorphism defined by $(x^1\ldots x^k)\to (-1)^{x_i}$.  We
define $D_g$ to be the sum of the $C_s$ such that
$\psi_s|_V=g$. Notice that the $D_g$ are disjoint and that $D:=\sum_gD_g$ is
the
union of
the $C_i$ that appear in $V$. If we write $v_j=(x^1_j\ldots x^k_j)$,  and
we
identify again $0,1\in \bbF_2$ with the integers $0,1$, then
it is not difficult to check that relations (\ref{fund}) can be rewritten
as:

\begin{equation}
2M_j\equiv\sum_i x_j^iC_i\ , \quad j=1\ldots r.
\end{equation}
So equations \ref{fund} can be solved uniquely by setting
$M_j=L_{v_j}$, $j=1\ldots r$. The corresponding cover
$\pi\colon Z\to Y$  satisfies conditions i) and ii) of the
statement. In addition, $Z$ is smooth by Proposition 3.1 of
\cite{ritaabel}, since $D$   is smooth, and it is connected since the set
of $g\in G$ such that $D_g\ne 0$ generates $G$. In order to complete the proof
we  have to check that for every
$v=(x_v^1\ldots x_v^k)\in V$ the eigensheaf   $M_v\inv$ of $\pi_*\OO_Z$ on
which
$G$ acts via the character
$v$ is $L_v\inv$. By Theorem $2.1$ of \cite{ritaabel}, we have
$2M_v\equiv\sum_g
\epsi_v(g)D_g$.  This equation can be rewritten as
$2M_v\equiv\sum_i x_v^iC_i$, and thus
$2L_v=2M_v$ in $\Pic(Y)$. The equality  $L_v=M_v$ follows since
$_2\Pic(Y)=0$.
\end{proof}

\begin{remark} Write $U:=Y\setminus
\cup_iC_i$\ .  Then there is an isomorphism $\psi\colon V \to
\, _2 \Pic(U)$ and the restriction to $U$ of
the cover $\pi\colon
Z\to Y$  is the $G-$torsor corresponding to $\psi$ under the natural map
$H^1(U,G)\to \Hom(V,\Pic(U))$.
\end{remark}
\begin{remark}  The proof of Proposition \ref{cover} shows that if one
removes the assumption
$_2\Pic(Y)=0$ then the cover $\pi\colon Z\to Y$ exists but it is not
determined
uniquely. Also, if one assumes
$char({\bf k})=2$, then the proof  shows the existence of a purely inseparable
cover with a $G-$action.
\end{remark}
\medskip

Let
$\eta\colon Y\to \Si$ be the map that contracts
  the curves $C_i$ that appear in $V$ to singular points of type $A_1$.
The inverse image in $Z$ of a curve $C_i$ that appears in $V$ is a
disjoint union of $2^{r-1}$
$(-1)-$curves. Blowing down all these $(-1)-$curves, we obtain a smooth
surface $\bar{Z}$ and
a $G-$cover $\bar{\pi}\colon \bar{Z}\to \Si$ \ branched precisely over the
singularities
of $\Si$. Then we have  the following commutative diagram: \[
\begin{CD}
Z@>  >>\bar{Z}\\
@V\pi VV@V\bar{\pi} VV\\
Y@>\eta>>\Sigma,
\end{CD}
\] We close this section by computing the invariants of $Z$ and
$\bar{Z}$.

\begin{lem}\label{c2} With the same assumptions and notations as in
Proposition
\ref{cover} (in
particular, $r$ is the dimension of $V$ and $m$ is the number of the
$C_i$ that
appear in $V$) one has:
$$c_2(Z)=2^rc_2(Y)-m2^r.$$
\end{lem}
\begin{proof}
If the base field is $\C$, then  the formula follows easily by topological
considerations. We give an algebraic proof, valid for fields of
characteristic
$\ne 2$.

Denote by $D$ the branch divisor of $\pi$ (which is the union of $m$
disjoint
nodal curves), and by $R=\pi\inv D$ the ramification divisor.  Consider the
following exact sequence of sheaves on
$Z$:
\begin{equation}\label{diff}
0\to\pi^*\Omega^1_Y\stackrel{j}{\to} \Omega^1_Z\to {\mathcal K}\to 0
\end{equation}
where the cokernel $\mathcal K$ is a torsion sheaf supported on $R$.
Consider a ramification  point $P\in Z$ and let $R'$ be the irreducible
component of $R$ containing $P$.  The subgroup $H\subset G$
consisting
of the elements that  induce the identity on $R'$ is isomorphic to $\Z_2$
(cf.
\cite{ritaabel}, Lemma $1.1$). The surface $W:=Z/H$ is smooth, since the
fixed
locus of $H$ is purely $1-$dimensional,  and
$\pi$ factorizes as
$Z\stackrel{\alpha}{\to}W\stackrel{\beta}{\to} Y$. Let $Q=\alpha(P)$ and
$D'=\alpha(R')$. The map $\beta$ is
\'etale in a neighbourhood of $Q$, and thus
$\beta^*\Omega^1_Y\hookrightarrow\Omega^1_W$ is an isomorphism locally near
$Q$. It follows that the inclusion
$\pi^*\Omega^1_Y\hookrightarrow\alpha^*\Omega^1_W$ is an isomorphism locally
around
$P$.     There exists an open
neighbourhood
$U$ of $Q$ in $W$ such that $Z|_U$ is defined in $U\times {\mathbb A}^1$ by
the
equation
$z^2=b$, where $b$ is a local equation for $D'$ and $z$ is the affine
coordinate in
${\mathbb A}^1$. Notice that $z$ is a local equation for $R'$ around $P$.
Let
$x$ be a function on
$W$ such that
$x,b$ are local parameters on $W$ around $Q$. Then  the map
$j$ of  sequence
(\ref{diff}) can be written locally  as $(dx,db)\mapsto(dx, 2zdz)$.  It
follows
that the cokernel $\mathcal K$ is naturally isomorphic to the conormal
sheaf of $R$, $\OO_R(-R)$. A standard computation with
Chern classes gives:
$c_2(Z)=2^rc_2(Y)+2R^2+\pi^*K_Y\cdot R=2^rc_2(Y)+2^{r-1}D^2+2^{r-1}K_Y\cdot
D=2^rc_2(Y)-m2^r$.

\end{proof}

\begin{prop}\label{invariants}
Under the same assumptions and notation as above the following holds:

 $$\kappa(Z)=\kappa(\bar{Z})=\kappa(Y);$$
$$K^2_Z=2^rK^2_Y-m2^{r-1} \quad K_{\bar{Z}}^2=2^rK^2_Y;$$
$$\chi(Z,\OO_Z)=\chi(\bar{Z},\calO_{\bar{Z}}) = 2^r\chi(\OO_Y)-m2^{r-3}.$$
\end{prop}

\begin{proof}
We have $K_{\bar{Z}}=\bar{\pi}^*K_{\Si}$, since $\bar{\pi}$ is
unramified in codimension
$1$ and $\Sigma$ is normal, and therefore
$K^2_{\bar{Z}}=2^rK^2_{\Si}=2^rK^2_Y$.  The formula for $K_Z^2$ follows
immediately.  Since $\chi$ is a birational invariant,  it is enough to
compute
it for $Z$. Then the
formula for $\chi(Z,\OO_Z)$ follows  from Lemma \ref{c2} and Noether's
formula.

If $\kappa(Z)=-\infty$, then we  have $\kappa(Y)=-\infty$ ($\pi$ is
separable since the characteristic of ${\bf k}$ is $\ne 2$). So assume that
$\kappa(Z)\ge 0$ and denote by
$\tilde{Z}$ the minimal model of
$Z$ and $\bar{Z}$. Then $G$ acts biregularly on $\tilde{Z}$. We denote by
$\tilde{\pi}\colon \tilde{Z}\to\tilde{\Si}:=\tilde{Z}/G$ the quotient map.
The
surface $\tilde{\Si}$ has canonical singularities and it is  birational to
$Y$ and
$\Si$. Denote by $\epsi\colon \tilde{Y}\to\tilde{\Si}$ the minimal
resolution.
We have
$K_{\tilde{Z}}=\tilde{\pi}^*K_{\tilde{\Si}}$ and thus $K_{\tilde{\Si}}$ and
$K_{\tilde{Y}}=\epsi^*K_{\tilde{\Si}}$ are  nef. So $\tilde{Y}$ is minimal
and, in addition, $K_{\tilde{Y}}\sim 0$ iff $K_{\tilde{Z}}\sim 0$ and
$K_{\tilde{Y}}^2=0$ iff $K_{\tilde{Z}}^2=0$. This remark shows that
$\kappa(\tilde{Y})=\kappa(\tilde{Z})$.
\end{proof}
\medskip

\section{Rational surfaces  with many nodes.}
Throughout this section we assume that $Y$ is a smooth rational
surface and $C_1,...,C_k$ are disjoint nodal curves of $Y$. As
before, we  let
$V$ be the code associated to the $C_i$, $r$ its dimension and $m$  the
number of the $C_i$
that
appear in $V$.
The group $\Pic(Y)$ is free abelian of rank $\rho(Y)=10-K^2_Y$ and the
intersection form on $\Pic(Y)$ induces a non degenerate $\bbF_2-$valued
bilinear
form on $\Pic(Y)/2\Pic(Y)$. Since $C_i^2=-2$ and the $C_i$ are disjoint,
the
image of
$\calC/2\calC$ is a totally isotropic subspace of $\Pic(Y)$. Thus the
dimension $r$ of $V$ satisfies $r\ge k-[\frac{\rho(Y)}{2}]$. As a corollary
of the results in the previous section we have the
following

\begin{lem}\label{emme}
If $r\ge 4$, then $m\ge 8$.
\end{lem}

\begin{proof} Consider the cover $\pi\colon
Z\to Y$  of Proposition \ref{cover} associated to $V$ and
the corresponding cover of $\Si$,
$\bar{\pi}\colon\bar{Z}\to \Si$. By Proposition
\ref{invariants},
$\bar{Z}$ is ruled and thus
$\chi(\bar{Z})\le 1$. The result follows by using the formula for
$\chi(Z,\OO_Z)$  of
Proposition \ref{invariants}.\end{proof}

\begin{thm}\label{doublefibres}
Let  $C_1\ldots C_k$ be disjoint nodal curves
on a rational surface
$Y$, let $V$ be the code associated to $C_1\ldots C_k$
 and assume that the length  of the reduced code $V'$  of
$V$ is
$m\ge 8$. Denote by $\eta\colon Y\to \Si$ the map that contracts to nodes
the
$C_i$ that appear in $V$. Then there exists a fibration
$\beta\colon
\Si\to
\pp^1$ such that:
\begin{itemize}
\item [(i)] the general fibre of $\beta$ is $\pp^1$;

\item [(ii)] $m=2n$ is even and $\beta$ has $n$ double fibres, each
containing two nodes   of
$\Si$;

\item [(iii)] the code $V$ is essentially isomorphic to $DE(n)$.
\end{itemize}
\end{thm}
\begin{proof} Let $\pi\colon
Z\to Y$ be the
cover of Proposition \ref{cover} and let
$\bar{\pi}\colon\bar{Z}\to \Si$\, be the
corresponding cover of $\Si$.  By Proposition
\ref{invariants}, the surface
$\bar{Z}$ is ruled
and  has irregularity
$q(\bar{Z})=1+m2^{r-3}-2^r>0$. Denote by $\alpha\colon \bar{Z}\to C$\, the
Albanese
pencil. By the canonicity of the Albanese map, the group $G$ preserves
the divisor class of a fibre. Consider the canonical
homomorphism $G\to \Aut(C)$: if it is not injective, then there exists
$g\in G$
that maps  each fibre of $\alpha$ to itself. Hence a general fibre, being
isomorphic to $\bbP^1$, has $2$ fixed points of $g$ and   the ramification
locus for the
action of $G$ has
components of dimension $1$, a contradiction since the $G-$cover is branched
precisely over the singularities of $\Si$.
 Thus we have a commutative diagram:
\[
\begin{CD}
\bar{Z}@>\bar{\pi}>>\Sigma\\
@V\alpha VV@V\beta VV\\
C@>p>>\bbP^1,
\end{CD}
\]
where $p\colon C\to\pp^1$ is a $G-$cover. The general fibre of
$\beta$ is
$\pp^1$, since it is isomorphic to the general fibre of $\alpha$.  Since the
genus of
$C$ is equal to
$q(\bar{Z})$, by the Hurwitz formula the branch locus of $p$ consists of
$n=m/2$ points
(the inverse image of a branch point consists of $2^{r-1}$ simple
ramification
points). The cover
$\bar{\pi}\colon
\bar{Z}\to
\Si$ is obtained from
$p$ by base change and normalization, thus the fibres  of $\beta$ over the
branch points
$y_1\ldots y_n$ of $p$ are of the form $f_i=2\delta_i$, $i=1\ldots n$,  and
$\cup
_i\delta_i$
contains  all the nodes of
$\Si$. We claim that each double fibre contains at least one node. Indeed,
otherwise $\delta_i$ would be contained in the smooth part of $\Si$ and so
it
would be a Cartier divisor with  $\delta_i^2=0$, $K_{\Si}\cdot \delta_i=-1$,
a
contradiction to the adjunction formula.

Set $\beta'=\beta\circ\eta$.
Then for every
$i$, one can write
${\beta'}^*y_i=2A_i+\sum_s C_{i, s}$ and it follows that for every choice
of $h\ne j$ the divisor $\sum_sC_{h,s}+\sum_tC_{j,t}$
is divisible by
$2$ in $\Pic(Y)$, namely it corresponds to a vector of $V$. Since the
weights of $V$
are all divisible by $4$, it follows easily that each $\delta_i$ contains
precisely $2$ nodes of $\Si$.  So it is possible to relabel the
$C_i$ in such a way that
$\beta'(C_{2j-1})=\beta'(C_{2j})=y_j$ for
$j=1\ldots n$,  and that $C_{2j}+C_{2j-1}+C_{2h}+C_{2h-1}$ is divisible by
$2$  in
$\Pic(Y)$ for every choice of $j,h$. This shows that
$V$ is essentially isomorphic to the code $DE(n)$.
\end{proof}

Next we apply the above results to describe rational surfaces
with ``many'' disjoint
nodal curves. We start by describing an example.

\begin{ex} Consider  a relatively minimal ruled rational surface 
$\F_e:= {\rm Proj}(\OO_{\pp^1}\oplus \OO_{\pp^1}(e))$, $e\ge 0$,  and a point
$y\in\F_e$. If one blows up $y$, then the total transform of the ruling of
$\F_e$
containing
$y$ is the union of two $(-1)-$curves $E$ and $E'$ that intersect
transversely in a point
$y_1$. If one blows  up also $y_1$, then the strict transforms of $E$ and
$E'$ are
disjoint
nodal curves. By repeating  this procedure $n$ times  at points lying on
different
rulings of $\F_e$, one obtains a rational surface $Y$ containing $2n$
disjoint nodal
curves. One has $\rho(Y)=2n+2$ and it is easy to check that the code $V$
associated to
this collection of curves is $DE(n)$. We will call $Y$ the
{\em standard example} of a
rational surface with $\rho(Y)-2$ disjoint nodal curves.
\end{ex}

\begin{thm}\label{MT} Let $Y$ be a smooth rational surface and let
$C_1\ldots C_k$ be disjoint nodal curves of $Y$. Then:
\begin{itemize}
\item [(i)] $k\le \rho(Y)-1$,  and equality holds if and only if
$Y=\bfF_2$;
\item [(ii)] if $k=\rho(Y)-2$ and $\rho(Y)\ge 5$, then $Y$ is the standard
example. In particular  $k=2n$ is even and  the code $V$ is
$DE(n)$.
\end{itemize}
\end{thm}
\begin{proof} The group $\Pic(Y)$ is  free abelian  of rank $\rho:=\rho(Y)$.
The intersection form on
$\Pic(Y)$ extends to  a nondegenerate bilinear form of signature $(1,
\rho(Y)-1)$
on
$N^1(Y):=\Pic(Y)\otimes\R$.
The subspace  of $N^1(Y)$ spanned by the classes of the $C_i$ has
dimension
$k$ and the intersection form is negative definite there, thus we get
$k<\rho$.

We start by proving  (ii).  As before, we let $m\le k$ be the number of
nodal curves
that appear in the code $V$.  Recall  that the dimension $r$ of $V$ is
$\ge
\rho-2-[\frac{\rho}{2}]= [\frac{\rho+1}{2}]-2$. So, for $\rho\ge 11$, we
have
$r\ge 4$ and
thus $m\ge 8$ by Lemma \ref{emme}.

Assuming then that
$\rho\ge 11$,  we can apply Theorem
\ref{doublefibres}. Thus $m$ is even, say $m=2n$,
and there exists a morphism $\beta\colon \Si\to\pp^1$ such that the
general fibre of
$\beta$ is $\pp^1$ and $\beta$ has $n$ double fibres,
occurring at points $y_1\ldots y_n$ of $\pp^1$. Each double fibre contains
precisely $2$ nodes of
$\Si$, and the code $V$ is $DE(n)$. So we have $n-1=\dim V\ge
[\frac{\rho+1}{2}]-2$, namely $\rho-2\ge m=2n\ge 2[\frac{\rho+1}{2}]-2$. It
follows that $\rho$ is even and
$\rho-2=m=2n$. In particular, $m=k$, i.e. all the $C_i$ appear in $V$.

Set $\beta'=\beta\circ \eta$, denote by $F$ the cohomology class on $Y$ of a
fibre
of
$\beta'$ and let
\[T = \{L\in N^1(Y):L\cdot F = 0\}.\]
 A
basis  of
$T$ is  given by
$F$ and the classes of
$C_1\ldots C_{2n}$, since these are independent classes and $\dim
T=\rho-1=2n+1$.
On the other hand,  it is well known that, if one removes a component from
each reducible
fibre of $\beta'$, then $F$ and the classes of the remaining components of
the reducible
fibres are independent. It follows that the $F_i:={\beta'}^* y_i$,
$i=1,\ldots n$ are the
only reducible fibres of $Y$. As in the proof of Theorem \ref{doublefibres},
it is
possible to relabel the $C_i$ in such a way that  for each
$i$ one has
$F_i=\lambda_iC_{2i-1}+\mu_iC_{2i}+2\nu_i D_i$\,, with $D_i$ irreducible and
such that $D_i^2<0$.  From
$K_Y\cdot F=-2$, we get
$\nu_i=1$, $K_Y\cdot D_i=-1$,   and thus
$D_i^2=-1$, namely $D_i$ is a $(-1)$-curve. The curve  $D_i$ has nonempty
intersection with both $C_{2i-1}$ and $C_{2i}$, since $F_i$ is connected. So
the
equality
$$0=D_i\cdot F=D_i\cdot (2D_i+\lambda_iC_{2i-1}+\mu_iC_{2i})=
-2+\lambda_iD_i\cdot C_{2i-1}+\mu_iD_i\cdot C_{2i}$$
gives:
$$\lambda_i=\mu_i=D_i\cdot C_{2i-1}=D_i\cdot C_{2i}=1.$$
Blowing down $D_1\ldots D_n$ one obtains a smooth
surface ruled over $\pp^1$  with precisely $n$ reducible rulings, each
consisting of two $(-1)$-curves intersecting transversely. Blowing
down a $(-1)-$curve of each ruling, we obtain a  ruled surface
${\bf F}_e$. So $Y$ is the standard example.

In order to complete the proof of (ii), we have to describe the cases
$5\le \rho\le 10$.
 In addition we may assume $m < 8$, since for $m=8$ (and $\rho=10$) one can
apply the argument above to  show that $Y$ is the standard example.
Since $m<8$  all the elements of $V$ have weight $4$ and it is easy to
check that the only (numerical) possibilities for
the pair
$(k,r)$ are:
$(4,1)$, $(6,2)$,  $(7,3)$ and $(8,3)$. One has $m=k$ in all cases but the
last one, where $m=7$.

 Consider  the first three cases. Let  $Z\to Y$ be the Galois cover
considered in Proposition \ref{cover}
 and $\bar{Z}\to \Si$ the corresponding cover of $\Si$.  By
Proposition
\ref{invariants},
$\bar{Z}$ is a surface  satisfying $\kappa(\bar{Z})=\kappa(Y)$,
$K^2_{\bar{Z}}=8$, $\chi(\bar{Z})=1$. So $\bar{Z}$ is rational and
$K^2_{\bar{Z}}=8$ implies that
$\bar{Z}=\F_e$ for some $e\ge 0$.
Denote by $t$ the
trace of $g\in G\setminus\{1\}$ on the $l$-adic cohomology
$H^2(\bar{Z},\bbQ_l) \cong \bbQ_l^2$. Since the class in $H^2(\bar{Z},\bbQ_l)$ of the
canonical bundle of $Y$ is $G-$invariant,  $t$  is either equal to $0$ or $2$. Applying
the ($l$-adic) Lefschetz fixed point formula (see \cite{sga},  (4.11.3), cf. the
next
section for the analogous statement for the complex cohomology) we see that
$t=0$ is impossible and hence $g$ acts identically on $H^2(\bar{Z},\bbQ_l)$.
In
particular, given the ruling (or a ruling if
$e=0$)  $f\colon \bar{Z}\to \pp^1$
 the action of the Galois group $G$ of
$\bar{Z}\to\Si$   descends to an action on
$\pp^1$ and  there is an induced fibration $h\colon
\Si\to \pp^1/G=\pp^1$. The same argument as in the proof  of Theorem
\ref{doublefibres}
shows that the action of
$G$ on
$\pp^1$ is faithful. Thus each element $g$
of
$G$ fixes precisely two fibres of $h$,
each containing two fixed points of $g$.  Since
$\text{Aut}(\bbP^1)$ does not contain a subgroup isomorphic to $\bbZ_2^3$, we
can rule out immediately
 the case $(n,r) = (7,3)$. In the remaining two cases, the cover
$\pp^1\to\pp^1/G$ is branched over
$n=m/2$ points,
and over each of these points $h$ has an irreducible double fibre
containing $2$ nodes. It follows easily that $Y$ is the standard
example.

Finally consider the case $(8,3)$, $m=7$ (the code is essentially
isomorphic to the Hamming code defined by the root lattice of type $E_7$).
By Proposition \ref{invariants} the
$G$-cover
$\bar{Z}\to
\Si$ is a smooth ruled surface with invariants $K^2_{\bar{Z}}=0$,
$\chi(\bar{Z})=1$. Thus $\bar{Z}$ is rational. The preimage of the
nodal curve of
 $Y$ not
appearing in $V$ is a set $D_1\ldots D_8$
of  disjoint
nodal curves  on which  $G=\Z_2^3$ acts transitively. The code
$\tilde V\subset \bbF_2^8$ associated to $D_1\ldots D_8$ is acted on by $G$,
and
therefore all the nodes appear in $\tilde V$, namely $\tilde V$ has  $m=8$.
Thus
$\bar{Z}$ is a standard example with $\rho=10$. If there is only one
pencil with rational fibres $f\colon \bar{Z}\to\pp^1$ such that the $D_i$
are
contracted by $f$, then one argues as in case $(7,3)$ and obtains a
contradiction by
showing  the existence of a
$\Z_2^3-$cover $\pp^1\to\pp^1$.
So assume that
there are two pencils with rational fibres
$f_j\colon \bar{Z}\to\pp^1$, $j=1,2$ such that the $D_i$ are contracted both
by
 $f_1$
and
$f_2$. Denote by
$F_j$,
$j=1,2$,   the class in $N^1(\bar{Z}):=\Pic(\bar{Z})\otimes\R$ of a smooth
fibre of
$f_j$. Considering
the
intersection form, one sees immediately that the classes of
$F_1$,
$F_2$, $D_1\ldots D_8$ are a basis of $N^1(\bar{Z})$. Consider a
nonzero
$g\in G$. The surface $Z':=\bar{Z}/<g>$ is a rational surface with $s$
singular points of type
$A_1$, that are the images of the fixed points of $g$ on $\bar{Z}$. By the
standard double cover formulas:
 $$1=\chi(\bar{Z})=2\chi(Z')-s/4=2-s/4$$ and so $s=4$.
Denote by $t$ the
trace of $g$ on the $l$-adic cohomology
$H^2(\bar{Z},\bbQ_l)$. Applying again the  Lefschetz fixed point
formula
 we get $t=2$.  The action of
$g$ on
$H^2(\bar{Z},\bbQ_l)$ preserves the subspace $<D_1\ldots D_8>$ generated by
the fundamental
classes of the divisors $D_1\ldots D_8$, and thus it preserves  also its
orthogonal subspace, which is
spanned by the classes of $F_1,F_2$. The trace of $g$ on $<D_1\ldots D_8>$
is zero. It follows
that
$g$ is the identity on $<F_1, F_2>$, namely every $g\in G$ preserves both
pencils. Thus we can apply  again the argument above to one of the pencils
and the
proof of (ii) is complete.

Finally we prove (i). Assume that $k=\rho(Y)-1$. The code $V$ has
length
$\rho-1$,  dimension
$r\ge
[\frac{\rho+1}{2}]-1$ and all the weights  divisible by $4$. Thus
if $\rho\ge
9$, then $m\ge 8$ by Lemma \ref{emme} and one can argue as in case (ii)
and show
that $Y$ is the surface constructed in the standard example and $V$ is
essentially
isomorphic to $DE(n)$, with $n=\rho/2-1$. In particular, $r=n-1=\rho/2-2$,
contradicting
$r\ge
[\frac{\rho+1}{2}]-1$. So assume $\rho\le 8$. If $\rho=2$, then
$K^2_Y=8$
and so $Y$ is the minimal ruled surface $\F_2$. If $\rho>2$,
the only numerical possibility is
$\rho=8$, $r=3$. Let $Y$ be a surface corresponding to this possibility. We
have $K^2_Y=2$. Up to a permutation, we may assume that
$C_1\ldots C_4$ is an even set. The corresponding double  cover
$Y'\to Y$ is a smooth rational surface (same proof as Proposition
\ref{invariants}), with $K^2_{Y'}=0$.  The inverse images of
$C_1\ldots C_4$ are $(-1)-$ curves, while the inverse images of
$C_5, C_6, C_7$ are three pairs of disjoint nodal curves.
Blowing  down the $(-1)-$curves, one obtains a rational surface
$Y''$ with $\rho(Y'')=6$ and containing $6$ disjoint nodal curves. This
is impossible, and the proof is complete.
\end{proof}

\begin{remark} For  a rational surface  $Y$ with $\rho(Y)\le 4$ containing
   $k=\rho(Y)-2$ disjoint nodal curves, the code $V$ is zero and one
 cannot
   argue as in Theorem 3.3. On the other hand, this case can be studied
   directly  and it is easy to check that the possibilities for
 $(k,\rho)$
   are:
\begin{itemize}
   \item[(i)]  $(0,2)$ and $Y$ is  a surface ${\bf F}_e$, $e\neq 2$.

  \item[(ii)]  $(1,3)$ and $Y$ the blowup of ${\bf F}_2$ at a point outside the
   negative section
  ( the nodal curve is the pull back of the negative section);
  or $Y$ is the blowup of ${\bf F}_1$ at a point on the negative section
  (the nodal curve is the strict transform of the negative
   section);

  \item[(iii)]  $(2,4)$ and $Y$ is
    the standard example with $k=2$; or $Y$ is the blowup of ${\bf F}_2$
at points
 $x_1$,
  $x_2$, with $x_1$ not
   on the negative section and $x_2$ infinitely near to $x_1$ ( the
   nodal curves  are the pullback of the negative section and the strict
  transform of the exceptional curve of the first blowup).
\end{itemize}

\end{remark}

\section{An application}
Throughout this section we  assume that the  ground field is
$\C$.
We  apply the previous results to study involutions (i.e.  automorphisms of
order 2)
on minimal surfaces of general type with $p_g = 0$ and $K_S^2 = 8$ or $9$.

We start by extending the results of section $3$ to complex surfaces
with $p_g=q=0$ and nonnegative Kodaira dimension. The use of  Miyaoka's
formula is a key ingredient for the proof below and explains the assumption that the ground field is
$\bbC$ in this section.

\begin{prop}\label{kappa>0}
Let $Y$ be a surface with $p_g(Y)=q(Y)=0$ and $\kappa(Y)\ge 0$, and let
$C_1\ldots
C_k\subset Y$ be disjoint nodal curves. Then:
\begin{itemize}
\item[(i)] $k\le \rho(Y)-2$;

\item[(ii)] if $k= \rho(Y)-2$, then $Y$ is minimal.
\end{itemize}
\end{prop}

\begin{proof} Assume first that $Y$ is minimal. In this case we
can apply Miyaoka's formula (\cite{miyaoka}, section $2$):
$3c_2(Y)-K^2_Y\ge \frac{9}{2}k$, and i) follows immediately using
$0\le K^2_Y\le 9$ and Noether's formula.

Now assume that $Y$ is not minimal and let $\bar{Y}$ be the minimal model of
$Y$.
We use induction on
$\nu:=\rho(Y)-\rho(\bar{Y})$.  Let $E\subset Y$ be an irreducible
$(-1)-$curve and let $Y'$ be the surface obtained by blowing down
$E$. If $E$ does not intersect any of the
$C_i$, then
$Y'$ contains
$k$ disjoint
nodal curves and induction gives: $k\le \rho(Y')-2=\rho(Y)-3$.
So assume, say, $C_1\cdot E=\alpha>0$. Then the image $C_1'$ of $C_1$ in
$Y'$ is
an irreducible curve
such that
$(C_1')^2=-2+\alpha^2$, $C_1'\cdot K_{Y'}=-\alpha$.
 Now necessarily
$\alpha=1$. In fact suppose that $\alpha\geq 2$. Then $C_{1}'{}^2>0$
and therefore  the image of $C_1'$ in the minimal model ${\bar
Y}$ of $Y$ is a curve $C_1''$. Since
$C_1''\cdot K_{\bar Y}\leq
C_1'\cdot K_{Y'}<0$ and $K_{\bar Y}$ is nef because  $\kappa(Y')\ge 0$, we
have
a contradiction. Therefore $C_1'$ is a $(-1)-$curve. In addition,
$E\cdot C_i=0$ for $i>1$, since otherwise
$Y'$ would contain a pair of irreducible $(-1)-$curves with
nonempty intersection, which is impossible again because
$\kappa(Y')\ge 0$. Now blowing down
$C_1'$ we obtain a surface $Y''$ containing a set of
$k-1$ disjoint irreducible nodal curves. Using induction again,
we have $k-1\le \rho(Y'')-2=\rho(Y)-4$ and the proof is complete.
\end{proof}

Let $S$ be a surface admitting an involution $\sigma$. Let $k$ be
the number of isolated fixed points of $\sigma$ and let $D$ be the
1-dimensional part of the fixed-point locus. The divisor $D$ is
 smooth  (possibly empty). If we consider
 the blow-up $X$ of the set of isolated fixed points, then the
involution
$\sigma$ lifts to an involution on $X$ (which we still denote by
$\sigma$) and  the quotient
$Y := X/<\sigma>$ has $k$ disjoint nodal curves $C_i$.\par We
recall
 the following
two well-known formulas:

(Holomorphic Fixed Point Formula) (see \cite{as}, pg.566):
\[\sum_{i=0}^2(-1)^i\text{Trace}(\sigma|H^i(S,\calO_S)) = \frac{k-D\cdot
K_S}{4}\]

(Topological Fixed Point Formula) (see \cite{gre}, (30.9)):
\[\sum_{i=0}^4(-1)^i\text{Trace}(\sigma|H^i(S,\bbC)) = k+e(D),\]
where $e(D) = -D^2-D\cdot K_S$ is the topological Euler
characteristic of $D$.
\begin{lem}\label{numeri}
Let $S$ be a surface with $p_g(S)=q(S)=0$ and let $\sigma$ be an
automorphism
of
$S$ of order $2$.
Let $D$ be the divisorial part of the fixed locus of $\sigma$, let $k$ be
the number of isolated
fixed points of $\sigma$ and let $t$ be the trace of $\sigma|H^2(S,\bbC)$.
Then:
$$k=K_S\cdot D+4;\quad t=2-D^2.$$

Furthermore  if $X$ is the blow-up of the $k$ isolated fixed
points of
$\sigma$, and $Y = X/<\sigma>$ one has $$\rho(S)+t = 2\rho(Y)-2k$$.
\end{lem}
\begin{proof}
The first fixed point formula gives
\begin{equation}\label{HFF}
k = 4+K_S\cdot D
\end{equation}
Together with the second formula we obtain
\begin{equation}\label{TFF}
t:=\text{Trace}(\sigma|H^2(S,\bbC)) = 2-D^2.
\end{equation}
For the last part notice that we have
$$e(S)+k = e(X) = 2e(Y)-2k-e(D).$$
Since by the topological fixed point formula $e(D) = -k+2+t$,
one has
$$e(S)+t+2  = 2e(Y)-2k.$$
Now $p_g=q=0$ implies $e(S) = \rho(S)+2, e(Y) = \rho(Y)+2$ and we obtain
\begin{equation}\label{B}
\rho(S)+t = 2\rho(Y)-2k
\end{equation}
\end{proof}

\begin{thm} A surface of general type $S$ with $p_g(S) = 0$ and $K^2_S =
9$ has no automorphism of order 2.
\end{thm}

\begin{proof} Assume otherwise. Since
$\rho(S) = 1$, we have
$t = 1$. Lemma \ref{numeri} gives  $D^2 = 1$. Since the canonical
class is invariant for $\sigma$, we have
$K_S \sim rD$ for some $r\in \bbQ$. Then $K_S^2=9$ yields  $K_S\sim 3D$
 and $K_S\cdot D = 3$.  Thus Lemma \ref{numeri} gives $k=7$ and
$2=2\rho(Y)-14$, i.e.
$\rho(Y) = 8$. So  $Y$ contains $\rho(Y)-1$ disjoint nodal
curves and $K_Y^2=2$. This is a contradiction in view of  Theorem
\ref{MT} and  Proposition \ref{kappa>0}.\end{proof}

\begin{thm}\label{class} Let $S$ be a minimal surface of general type with
$p_g(S) = 0$, $K_S^2 = 8$ and
let
$\sigma$ be an automorphism  of $S$ of order $2$.
Let $D$ be the divisorial part of the fixed locus of $\sigma$, let $k$ be
the number of isolated
fixed points of
$\sigma$ and  let
$Y$ be a minimal resolution of the quotient $S/<\sigma>$.
Then: $$D^2=0, \quad K_S\cdot D=k-4$$  and
 one of the following cases occurs:
\begin{itemize}
\item [(i)] $k = 4$, $D=0$ and   $Y$ is a minimal surface
of general type with $p_g(Y) = 0$ and $K_Y^2 = 4$.
\item [(ii)]   $k=6$,  and   $Y$ is a minimal surface of
general type with $p_g(Y) = 0$ and $K_Y^2 = 2$.
\item[(iii)] $k = 8$, $Y$ is a  minimal   surface with
$p_g(Y) = q(Y)= 0$, $\kappa(Y)=1$ for which the elliptic fibration $Y\to
\pp^1$ has two reducible fibres of Kodaira type $I_0^*$, and as such constant moduli.
\item [(iv)]$k = 10$,  and  $Y$ is a rational
surface  from Example 1 with $\rho =12$. The fibration with
connected rational fibres  $f\colon Y\to \pp^1$ pulls back on $S$ to a
pencil of  hyperelliptic  curves of genus $5$.
\item [(v)]$k = 12$,  and  $Y$ is a rational
surface  from Example 1 with $\rho =14$. The fibration with
connected rational fibres  $f\colon Y\to \pp^1$ pulls back on $S$ to a
pencil of  hyperelliptic curves of genus $3$.
\end{itemize}
\end{thm}

\begin{proof} Since $\rho(S) =
2$, the possible values for the trace $t$ are $0$ and
$2$. \par The case  $t = 0$ does not occur. Indeed, assume otherwise.  By
Lemma \ref{numeri}, $D^2 = 2$ so that $D\ne 0$. Since $t=0$,  the
invariant part of
$H^2(S,\bbQ)$ is one-dimensional and thus (because the canonical
class is invariant for $\sigma$),
$K_S \sim rD$ for some $r\in \bbQ$.  Thus $K_S \sim 2D$ and, hence
$K_S\cdot D = 4$. Lemma \ref{numeri} gives $k = 8$ and
$\rho(Y) = 9$, and so by Noether's formula $K_Y^2 = 1$. Since $Y$ contains $
8$
disjoint
 nodal curves, we have a contradiction to  Theorem
\ref{MT} and Proposition \ref{kappa>0}. So $t\neq 0$.

Now we consider the case  $t = 2$,
that is, the involution
$\sigma$ acts identically on $H^2(S,\bbQ)$. In this case $D^2 =
0$.

If $D= 0$, we get $k = 4$ and the surface $Y$ is a
surface of general type with
$K_Y^2 = 4$ and $\rho(Y) = 6$. It contains an even set of four
disjoint nodal curves $C_1,\ldots C_4$ and thus  it is minimal by
Proposition
\ref{kappa>0}. This is case (i).

The last case to consider is $t = 2$ and $D\ne 0$. Since
$D^2=0$, we have
$K_S\cdot D = 2m$, with $m>0$. Then Lemma \ref{numeri} gives $k = 4+2m$, so
that
in particular   $k$ is  $\ge 6$ and even, and
$\rho(Y) = 6+2m = k+2$.

Assume that $\kappa(Y)\ge 0$.
Since $Y$ is a minimal surface by Proposition
\ref{kappa>0}, $K_Y^2\geq 0$ and so
$k =8-K^2_Y\le 8$. So either $k=6$ or $k=8$. If $k=6$,
$K^2_Y=2$ and so
$Y$ is of general type and we have  case (ii). If $k=8$, then $K_Y^2=0$ and
thus
$Y$, being minimal,
is not of general type.
 Since $p_g(Y) = q(Y)= 0$, $Y$ is  either an Enriques surface or a surface
of
Kodaira dimension 1. The first case cannot occur. In fact since
$K_Y\sim 0$ and $D^2=0$ we would have $K_S\cdot D=0$, a
contradiction. So $\kappa(Y)=1$ and $Y$ is a minimal properly
elliptic surface.
 Denote by $f\colon Y\to \pp^1$ the elliptic fibration
and let $F$ be a general fibre of $f$. Since $K_Y$ is numerically
a rational multiple of
$F$, we have $F\cdot C_i=0$ for every $i$, namely the $C_i$ are mapped to
points
by $f$. Let
$\bar{F}$ be a fibre containing,  say,    $C_1\ldots C_s$ and let $A_1\ldots
A_p$ be the
remaining irreducible components of $\bar{F}$. It is well known that the
classes of
$A_1\ldots A_p, C_1\ldots C_r$ in $H^2(Y,\Q)$ are independent and span a
subspace $U_1$ on which
the intersection form is seminegative. The classes of $C_{s+1}\ldots C_{8}$
are also
independent  and span a subspace $U_2$ such that the intersection form is
negative on $U_2$ and
$U_1\cap U_2=\{0\}$. Since $\rho(Y)=10$, we see that the only possibility is
$p=1$. Looking at
Kodaira's list of singular elliptic fibres (see e.g. \cite{bpv}, pg.150),
one
sees that the possible types of singular
 fibres containing some of the $C_i$ are $I_2$, $I_0^*$ and $III$.
In addition, we have $12=e(Y)=\sum_te(F_t)$, where $F_t$ is the fibre of $f$
over
the point $t\in\pp^1$ and $e$ denotes the topological Euler--Poincar\'e
characteristic. It is
easy to check that the only numerical possibility is that $f$ has two
$I_0^*$ fibres, each
containing
$4$ of the
$C_i$, and that every other singular fibre is a multiple of a smooth
elliptic curve.
Up to a permutation we may assume that the $I_0^*$ fibres of $f$ are
$C_1+\ldots+C_4+2D_1$ and
$C_5+\ldots +C_8+2D_2$. So $C_1+ \ldots +C_8\equiv 2(F-
D_1-D_2)$ is divisible by $2$ in $\Pic(Y)$.
Let $\pi\colon Y'\to Y$ be the corresponding double cover. For a general
fibre $F$ of $f$,
$\pi^*F$ is disconnected and the Stein factorization of $f\circ \pi$ gives
rise to an
elliptic fibration $f'\colon  Y'\to \pp^1$ ``with the same fibres'' as $f$.
The inverse images of $D_1$, $D_2$ are smooth elliptic curves. The inverse
images of
$C_1,\ldots C_8$ are
$8$
$(-1)-$ curves contained in the fibres of
$f'$.
Blowing these exceptional curves down, one obtains an elliptic fibration
$f''\colon
Y''\to\pp^1$ whose only singular fibres are multiples of smooth elliptic
fibres. Thus $f''$
has constant moduli, and therefore $f'$ and   $f$ have constant moduli too.
This is case (iii).

Finally, assume that $Y$ is a rational surface. Since $k\geq 6$ and
$\rho(Y)=k+2$  we can
apply Theorem \ref{MT}
to obtain that $Y$ is as in the standard example. In particular there is a
fibration $f\colon
Y\to\pp^1$ with general fibre $F$ isomorphic to $\pp^1$.  If we write
$K_S\cdot D=2m$ (hence
$k=2m+4$), then $f$ has precisely $m+2$ singular fibres of the form
$C_{2i-1}+C_{2i}+2E_i$,
with $E_i$ a
$(-1)-$curve and $E_i\cdot C_{2i-1}=E_i\cdot C_{2i}=1$. Denote by
$\bar{D}$ the image of
$D$ on
$Y$ and by
$L$ the line bundle of $Y$ such
that $2L\equiv \bar{D}+C_1+\ldots +C_{k}$. The intersection number
$E_i\cdot \bar{D}=E_i\cdot 2L-E_i\cdot (C_1+\ldots C_{k})=2L\cdot E_i-2$
is even. Thus we may
write
$\bar{D}\cdot F=\bar{D}\cdot (2E_i+C_{2i-1}+C_{2i})=2\bar{D}\cdot E_i=4d$,
 and the pre-image in $X$ of the ruling
on $Y$  is a pencil of hyperelliptic curves  of genus $2d-1$.
Blowing down the curves $E_i$ and then the images of the $C_{2i}$, we obtain
a
birational morphism
$p\colon Y\to
\F_e$ onto a relatively minimal ruled
  surface.  Let $C$ be the image of $\bar{D}$ on $\F_e$. Let
$F,S$ be the standard generators of $\Pic(\F_e)$ with $F^2 = 0,
S^2 = -e \le 0, F\cdot S = 1$. We have $C \sim aF+4dS$. The curve
$C$ has $m+2$ singular points of type $(2d,2d)$, that are solved by the
morphism $p$.
 Since $\bar{D}^2=D^2/2 = 0$, we get
\[0=C^2- (m+2)8d^2 = 8d(a-2de-d(m+2)).\]
This gives us a first equation:
\begin{equation}
a = d(m+2)+2de.
\end{equation}
We also know that $\bar{D}\cdot K_Y = 2m$. On the other hand,
\[\bar{D}\cdot K_Y = C\cdot K_{\F_e}+2(m+2)2d,\]
and we get the second equation
\begin{equation}
a = 2d(m+e)-m.
\end{equation}
Comparing the two equations, we get $dm = m+2d$. This has the solutions
$(m,d) =  (3,3), (4,2)$, which yield the cases (iv) and (v), respectively.
\end{proof}

\bigskip
 \begin{remark} We do not know whether all the possibilities in Theorem
\ref{class} really
occur. One can check that in the case of the
bicanonical involution of the surface $S$ of  example (4.2) of
\cite{mp} the quotient is as in case (v). In addition,
$\text{Aut}(S)=\Z_2^3$ and the
remaining  involutions are as in case (iii). Example (4.3) of \cite{mp} has
a group $\Ga$ of  automorphisms isomorphic to $\Z_2^4$: some elements of
$\Ga$ have no
$1-$dimensional fixed part, and thus are as in case (i), while the others
are as
in case (iii).  Both examples are Beauville-type surfaces (cf
\cite{bpv}, pg. 236). We intend to return to this problem in a future
paper.
\end{remark}

\bigskip
\bigskip

\begin{minipage}{12.5cm}
\parbox[t]{5.5cm}{Igor Dolgachev\\
Department of Mathematics\\
University of Michigan\\
Ann Arbor, MI48109\\
 USA\\
idolga@math.lsa.umich.edu}
\hfill
\parbox[t]{5.8cm}{Margarida Mendes Lopes\\
CMAF\\
Universidade de Lisboa\\
Av. Prof. Gama Pinto, 2\\
1649--003 Lisboa, PORTUGAL\\
mmlopes@lmc.fc.ul.pt}\\
\bigskip
\bigskip

\parbox{5.5cm}{Rita Pardini\\
Dipartimento di Matematica\\
Universit\`a di Pisa\\
Via Buonarroti, 2\\
56127 Pisa, ITALY\\
pardini@dm.unipi.it}
\end{minipage}

\end{document}